\documentclass[12pt,
]{article}

\usepackage{amssymb,amsmath,amstext,amsthm}
\usepackage{mathbbol}
\usepackage{mathrsfs}

\newcommand{\nek}{\newcommand}
\nek{\renek}{\renewcommand}
\nek{\vyk}[1]{}

\DeclareFontFamily{U}{ltt}{}
\DeclareFontShape{U}{ltt}{m}{n}{%
  <5><6><7><8><9> 
   <->[1.095] gen*cmtt%
  <10><10.95><12><14.4><17.28><20.74><24.88> 
   <->[1.095] cmtt10}{}

\makeatletter
\if@twoside%
\else\fi%
\makeatother

\nek{\sekd}{\subsection}

\nek{\imar}[1]{\marginpar
[%
\vspace{0ex}%
\flushright\scriptsize\sl
{$\mtho\to$}\\
\vspace{-0.8ex}
{#1}%
]
{%
\vspace{0ex}%
\flushleft\scriptsize\sl
{$\mtho\longleftarrow$}\\
\vspace{-0.8ex}
{#1}%
}}
\nek{\itsep}{\itemsep=0.3ex plus 0.15ex minus 0.1ex}
\nek{\tenu}[1]{
\itsep}
\nek{\enrm} 
{\itsep}
\nek{\enar} {\tenu{{\rm(\arabic{enumi})}}}

\newtheorem{theore}            {Theorem}  
\newtheorem{propo}     [theore]{Proposition}
\newtheorem{lemm}      [theore]{Lemma}

\theoremstyle{definition}
\newtheorem*{prF}{Proof}

\nek{\bpro}{\begin{propo}}
\nek{\epro}{\end{propo}}
\nek{\ble} {\begin{lemm}}
\nek{\ele} {\end{lemm}}
\nek{\bte} {\begin{theore}}
\nek{\ete} {\end{theore}}
\nek{\bpf} {\begin{prF}}
\nek{\epf} {\qed\end{prF}}
\nek{\qeD} [1] {\qed\;({\sl#1\/})} 
\nek{\epF} [1] {\qeD{#1}\end{prF}}

\nek{\ben}{\begin{enumerate}\itsep}
\nek{\een}{\end{enumerate}}
\nek{\bay}{\begin{array}}
\nek{\eay}{\end{array}}

\nek{\ltt} {\usefont{U}{ltt}{m}{n}}

\nek{\ie}{i.\,e.}

\nek{\tsup} {\mathop{\text{\ltt sup}}}
\nek{\ran} {\mathop{\text{\ltt ran}}}
\nek{\Ord}  {{\text{\ltt Ord}}}

\nek{\ZFC}  {{\text{\bf ZFC}}}
\nek{\ZF}   {{\text{\bf ZF}}}

\nek{\OD} {{\text{\rm OD}}} 
\nek{\ROD}{{\text{\rm ROD}}} 

\nek{\al} {\alpha}
\nek{\ba} {\beta}
\nek{\ga} {\gamma}
\nek{\ka} {\kappa}
\nek{\la} {\lambda}
\nek{\vpi}{\varphi}
\nek{\om} {\omega}
\nek{\omi} {\om_1}
\nek{\gc} {{\mathfrak c}}

\newlength{\spacel}
\nek{\dd}[1]
{\settowidth{\spacel}{\ }%
{\mtho$\hspace{0.2ex}{#1}$}-\linebreak[0]\hspace*{-1\spacel}\ }

\nek{\BBB}{} 
\nek{\dN}{{\BBB{\mathbb N}\BBB}}
\nek{\dR}{{\BBB{\mathbb R}\BBB}}
\nek{\skr}{\mathscr}
\nek{\cL} {{\skr L}}
\nek{\cM}  {{\skr M}}
\nek{\cP} {{\skr P}}
\nek{\cR}  {{\skr R}}
\nek{\cX} {{\skr X}}

\nek{\res} {\mathop{\restriction}}
\nek{\dm}  {$$}
\nek{\sus} {\mathopen{\exists\hspace{0.35ex}}}
\nek{\kaz} {\mathopen{\forall\hspace{0.35ex}}}
\nek{\imp} {\Longrightarrow}
\nek{\mpi} {\Longleftarrow}
\nek{\eqv} {\Longleftrightarrow}
\nek{\mo}  {\models}
\nek{\sq}  {\subseteq}
\nek{\nin} {\not\in}
\nek{\limp}{\,\imp\,}
\nek{\leqv}{\,\eqv\,}
\nek{\mto} {\longmapsto}
\nek{\ang} [1] {\langle #1\rangle}
\nek{\stk} [2] {{\ang{#1\hspace{0.3ex};\hspace{0.1ex}#2}}}
\nek{\ans} [1] {\{\hspace{0.01ex}#1\hspace{0.01ex}\}}
\nek{\itla} {\item\label}

\nek{\zi} {,\linebreak[0]\,}
\nek{\zd} {,\linebreak[0]\:}
\nek{\zt} {,\linebreak[0]\;}

\newlength{\dxii}
\nek{\mtho}{\mathsurround=0mm}%
\nek{\msur}{\hspace{-1\mathsurround}}
\nek{\hsur}{\hspace{-0.5\mathsurround}}
\nek{\noi}{\noindent}
\nek{\vom}{\vspace{1mm plus 0.1mm minus 0.3mm}}
\nek{\vtm}{\vspace{2mm plus 0.2mm minus 0.6mm}}

\nek{\dagg} {{\mtho(\ensuremath\dag)}}
\nek{\dagd} {{\mtho(\ensuremath\ddag)}}
\nek{\mast} {{\mtho(\ensuremath\ast)}}

\nek{\supp} [1] {||#1||}

\nek{\adR}  {{{\vphantom{X^x}}^\ast\hspace{-0.25ex}\dR}}
\nek{\ar}   {{^\ast\hspace{-0.25ex}r}}
\nek{\ax}   {{^\ast\hspace{-0.25ex}x}}
\nek{\aq}   {{^\ast\hspace{-0.4ex}q}}

\nek{\fy} {\mathbf y}
\nek{\fx} {\mathbf x}

\nek{\dwa}{\downarrow}
\nek{\upa}{\uparrow} 
\nek{\xle}{<_{\text{\tt lex}}}

\nek{\qU}{\mathopen{U}\hspace{0.1ex}}
\nek{\qD}{\mathopen{D}\hspace{0.1ex}}

\nek{\kla}[1] {[#1]_D}

\nek{\los} {\L o\v s}

\nek{\aE} {{^\ast\hspace*{-0.45ex} E}}

\nek{\eD} {E^D}
\nek{\eqD} {=^D}

\nek{\aR} {{{\vphantom{X^x}}^\ast\hspace*{-0.6ex}\cR}}
\nek{\aM} {{{\vphantom{X^x}}^\ast\hspace*{-0.6ex} M}}
\nek{\acM} 
{{{\hspace*{0.5ex}\vphantom{X^x}}^\ast\hspace*{-1.65ex}\cM}}

\begin{document}

\title{A definable nonstandard model of the reals}

\author{Vladimir Kanovei\thanks
{Partial support of RFFI grant 03-01-00757 and DFG 
grant acknowledged.} 
\and 
Saharon Shelah\thanks
{Supported by The Israel Science Foundation. Publication 825.}
}

\date{August 2003}
\maketitle

\begin{abstract}
We prove, in $\ZFC,$ the existence of a definable,  
countably saturated elementary extension of the reals. 
\end{abstract}

\sekd*{Introduction}

It seems that it has been taken for granted that there is 
no distinguished, definable nonstandard model of the reals.
(This means a countably saturated elementary extension of 
the reals.)
Of course if ${\bf V }  = {\bf L}$ then there is such an 
extension  
(just take the first one in the sense of the canonical 
well-ordering of $\bf L$), 
but we mean the existence provably in $\ZFC.$ 
There were good reasons for this: 
without Choice we cannot prove the existence of {\it any\/} 
elementary extension of the reals containing an infinitely 
large integer.~\footnote 
{In fact, from any nonstandard integer we can define a  
non-principal ultrafilter on $\dN,$ even a Lebesgue 
non-measurable set of reals \cite{lux}, 
yet it is consistent with $\ZF$ 
(even plus Dependent Choices) 
that there are no such ultrafilters as well as 
non-measurable subsets of $\dR$ \cite{sol}.}
\footnote
{It is worth to be mentioned that definable nonstandard 
elementary extensions of $\dN$ do exist in $\ZF.$ 
For instance, such a model can be obtained in the form of 
the ultrapower $F/U,$ where $F$ is the set of all 
arithmetically definable functions $f:\dN\to\dN$ while $U$ 
is a non-principal ultrafilter in the algebra $A$ of 
all arithmetically definable sets $X\sq\dN$.} 
Still there is one. 
                
\bte[{{\rm$\ZFC$}}]
\label t
There exists a definable,  
countably saturated extension  $\adR$ of the 
reals\/ $\dR,$ elementary 
in the sense of the language containing a symbol for every 
finitary relation on\/ $\dR$.
\ete

The problem of the existence of a definable proper elementary 
extension of $\dR$ was communicated to one of the authors 
(Kanovei) by V.\,A.\,Uspensky. 

A somewhat different, but related problem of 
{\it unique existence\/} of a  nonstandard real line $\adR$ 
has been widely discussed by specialists in nonstandard 
analysis.~\footnote
{``What is needed is an underlying set theory which proves the 
unique existence of the hyperreal number system [\dots]'' 
(Keisler~\cite[p.~229]k). 
} 
Keisler notes in \cite[\S\,11]k that, for any cardinal $\ka,$ 
either inaccesible or satisfying $2^\ka=\ka^+,$ there exists 
unique, up to isomorphism, \dd\ka saturated nonstandard real 
line $\adR$ of cardinality $\ka,$ 
which means that a reasonable level of 
uniqueness modulo isomorphism can be achieved, say, under 
GCH. 
Theorem~\ref t provides a countably saturated 
nonstandard real line $\adR,$ unique in absolute sense by virtue 
of a concrete definable construction in $\ZFC.$  
A certain modification of this example also admits a 
reasonable model-theoretic 
characterization up to isomorphism (see Section~\ref{s4}). 

The proof of Theorem~\ref t is a combination of several 
known arguments. 
First of all (and this is the key idea), 
arrange all non-principal ultrafilters 
over $\dN$ in a linear order $A,$ 
where each ultrafilter appears repetitiously as $D_a\zt a\in A.$ 
Although $A$ is not a well-ordering, we can apply the iterated 
ultrapower construction in the sense of \cite[6.5]{ck} 
(which is ``a finite support iteration'' in the forcing 
nomenclature), to obtain an ultrafilter $D$ in the algebra of 
all sets $X\sq \dN^A$ concentrated on a finite number of axes 
$\dN.$ 
To define a \dd Dultrapower of $\dR,$ the set $F$ of all 
functions $f:\dN^A\to\dR,$ also concentrated on a finite number 
of axes $\dN,$ is considered.   
The ultrapower $F/D$ is $\OD,$ 
thar is, ordinal-definable, actually, definable 
by an explicit construction in $\ZFC,$ hence, we obtain an  
$\OD$ proper elementary extension of $\dR.$ 
Iterating the \dd Dultrapower construction $\omi$ times  
in a more ordinary manner, \ie, with direct limits 
at limit steps, we obtain a definable 
countably saturated extension. 

To make the exposition self-contained and available for a 
reader with only fragmentary knowledge of ultrapowers, we 
reproduce several well-known arguments instead of 
giving references to manuals.

\sekd{The ultrafilter}
\label{s1}
                          
As usual, $\gc$ is the cardinality of the continuum. 

Ultrafilters on $\dN$ hardly admit any definable linear 
ordering, but maps $a:\gc\to\cP(\dN),$ whose ranges are 
ultrafilters, readily do. 
Let $A$ consist of all maps $a:\gc\to\cP(\dN)$ such that 
the set $D_a=\ran a=\ans{a(\xi):\xi<\gc}$ is an 
ultrafilter on $\dN.$ 
The set $A$ is ordered lexicographically: 
$a\xle b$ means that there exists $\xi<\gc$ such that 
${a\res\xi}={b\res\xi}$ and $a(\xi)<b(\xi)$ in the sense 
of the lexicographical linear order $<$ on $\cP(\dN)$ 
(in the sense of the identification of any $u\sq\dN$ 
with its characterictic function). 

For any set $u,$ $\dN^u$ denotes the set of all maps 
$f:u\to\dN$. 

Suppose that $u\sq v\sq A.$ 

If $X\sq \dN^v$ then put 
$X\dwa u=\ans{x\res u:x\in X}$. 

If $Y\sq \dN^u$ then put 
$Y\upa v=\ans{x\in \dN^v:x\res u\in Y}$. 

We say that a set $X\sq \dN^A$ is {\it concentrated\/} on 
$u\sq A,$ if $X=(X\dwa u)\upa A;$ in 
other words, this means the following:
\dm
\kaz x \zi y \in \dN^A\;
\big(
{x \res u= y \res u}\limp{({x \in X}\eqv{ y \in X})}
\big)\,.
\eqno(\ast)
\dm 
We say that $X$ is a {\em set of finite support\/}, if it is 
concentrated on a finite set $u\sq A.$  
The collection $\cX$ of all sets $X\sq\dN^A$ of finite 
support is closed under unions, intersections, complements, 
and differences, \ie, it is an algebra of subsets of $\dN^A.$
Note that if \mast\ holds for finite sets $u\zi v\sq A$ then 
it also holds for $u\cap v.$ 
(If $x \res{(u\cap v)}= y \res{(u\cap v)}$ then consider 
$ z \in \dN^A$ such that $ z \res u=x \res u$ and 
$ z \res v= y \res v$.) 
It follows that for any $X\in\cX$ there is a least 
finite $u=\supp X\sq A$ satisfying \mast.

In the remainder, 
if $U$ is any subset of $\cP(I),$ where $I$ is a given set, 
then $\qU i\:\Phi(i)$ 
({\it generalized quantifier\/})
means that the set $\ans{i\in I:\Phi(i)}$ belongs to $U.$ 

The following definition realizes the idea of a finite 
iteration of ultrafilters. 
Suppose that $u={a_1<\dots<a_n}\sq A$ is a finite set. 
We put
\dm
\bay{lcl}
D_u &=&\ans{X\sq \dN^u: D_{a_n}k_n\dots D_{a_2}k_2 \:D_{a_1}k_1\:
(\ang{k_1,k_2,...,k_n}\in X)}\,;\\[1.3\dxii]

D &=&\ans{X\in\cX:X\dwa{\supp X}\in D_{\supp X}}
\,.
\eay
\dm
The following is quite clear. 

\bpro
\label1
\ben
\enrm
\itla{11}
$D_u$ is an ultrafilter on\/ $\dN^u\;;$ 

\itla{12}
if\/ $u\sq v\sq A,$ $v$ finite, $X\sq\dN^u,$ then\/ 
$X\in D_u$ iff\/ $X\upa v\in D_v\;;$

\itla{13}
$D\sq\cX$ is an ultrafilter in the algebra\/ $\cX\;;$ 

\itla{14}
if\/ $X\in\cX,$ $u\sq A$ finite, and\/ $\supp X\sq u,$ then\/ 
${X\in D}\leqv{X\dwa u\in D_u}$.\qed
\een
\epro

\sekd{The ultrapower}
\label{s2}

To match the nature of the algebra $\cX$ of sets $X\sq\dN^A$ 
of finite support, we consider the family $F$ of all 
$f:\dN^A\to \dR,$ {\it concentrated\/} on some finite set $u\sq A,$ 
in the sense that  
\dm
\kaz x \zi y \in \dN^A\;
\big(
{x \res u= y \res u}\limp{f(x )=f( y )}
\big)\,.
\eqno\dagg
\dm 
As above, for any $f\in F$ there exists a least finite
$u=\supp f\sq A$ satisfying \dagg. 

Let $\cR$ be the set of all finitary relations on $\dR.$ 
For any \dd nary relation $E\in\cR$ and any 
$f_1,...,f_n\in F,$ define                              
\dm
\eD(f_1,...,f_n)\,\leqv\,
\qD\,x \in \dN^A\:E(f_1(x ),...,f_n(x ))
\,.
\dm
The set 
$X=\ans{x \in \dN^A:E(f_1(x ),...,f_n(x ))}$ 
is obviously concentrated on $u=\supp{f_1}\cup\dots\cup\supp{f_n},$ 
hence, it belongs to $\cX,$ and 
$\supp X\sq u=\supp{f_1}\cup\dots\cup\supp{f_n}.$ 

In particular, $f\eqD g$ means that 
$\qD\,x \in \dN^A\:(f(x )=g(x )).$ 
The following is clear:

\bpro
\label3
$\eqD$ is an equivalence relation on\/ $F,$ 
and any relation on\/ $F$ of the form\/ $\eD$ 
is\/ \dd{\eqD}invariant.\qed
\epro

Put $\kla f=\ans{g\in F:f=^D g},$ 
and $\adR=F/D=\ans{\kla f:f\in F}.$ 
For any \dd nary ($n\ge1$) relation $E\in\cR,$ let 
$\aE$ be the relation on $\adR$ defined as follows:
\dm
\aE(\kla{f_1},...,\kla{f_n})
\quad\text{iff}\quad
\eD(f_1,...,f_n)
\quad\text{iff}\quad
\qD\,x \in \dN^A\:E(f_1(x ),...,f_n(x )). 
\dm
The independence on the choice of representatives in the 
classes $\kla{f_i}$ follows from Proposition~\ref3.
Put $\aR=\ans{\aE:E\in\cR}.$ 
Finally, for any $r\in\dR$ we put $\ar=\kla{c_r},$ where 
$c_r\in F$ satisfies $c_r(x )=r \zd\kaz x $.

Let $\cL$ be the first-order language containing a symbol $E$ 
for any relation $E\in\cR.$ 
Then $\stk\dR\cR$ and $\stk\adR\aR$ are \dd\cL structures. 

\bte
\label e
The map\/ $r\mto\ar$ is an elementary embedding 
(in the sense of the language\/ $\cL$) of the 
structure\/ $\stk\dR\cR$ into\/ $\stk\adR\aR$. 
\ete
\bpf
This is a routine modification of the ordinary argument. 
By $\cL[F]$ we denote the extension of $\cL$ by functions 
$f\in F$ used as parameters. 
It does not have a direct semantics, but if $\vpi$ is a 
formula of $\cL[F]$ and $x \in\dN^A$ then $\vpi[x ]$
will denote the formula obtained by the substitution 
of $f(x )$ for any $f\in F$ which occurs in $\vpi.$ 
Thus,
$\vpi[x ]$ is an \dd\cL formula with parameters in $\dR$.

\ble[{{\rm\los}}]
\label l
For any closed\/ \dd{\cL[F]}formula\/ $\vpi(f_1,...,f_n)$ 
(all parameters\/ $f_i\in F$ indicated), we have$:$ \ 
\dm
\stk\adR\aR\mo\vpi(\kla{f_1},...,\kla{f_n})\;\leqv\; 
\qD\,x \;(\stk\dR\cR\mo\vpi(f_1,...,f_n)[x ]).
\dm 
\ele
\bpf
We argue by induction on the logic complexity of $\vpi.$ 
For $\vpi$ an atomic relation $E(f_1,...,f_n),$ the result 
follows by the definition of $\aE.$ 
The only notable induction step is $\sus$ in the direction 
$\mpi.$ 
Suppose that $\vpi$ is $\sus y\:\psi(y,f_1,...,f_n),$ 
and 
\dm
\qD\,x \;(\stk\dR\cR\mo\vpi(f_1,...,f_n)[x ]),
\quad\text{that is,}\quad
\qD\,x \;(\stk\dR\cR\mo\sus y\:\psi(y,f_1,...,f_n)[x ])\,.
\dm
Obviously there exists a function $f\in F,$ concentrated on 
$u=\supp{f_1}\cup\dots\cup\supp{f_n},$ such that, for any 
$x \in\dN^A,$ if there exists a real $y$ satisfying 
$\stk\dR\cR\mo\psi(y,f_1,...,f_n)[x],$ then $y=f(x)$ also 
satisfies this formula, \ie,  
$\stk\dR\cR\mo\psi(f,f_1,...,f_n)[x].$
Formally, 
\dm
\kaz x \in \dN^A\;
\big(
\sus y\in\dR\;(\stk\dR\cR\mo\psi(y,f_1,...,f_n)[x ])
\;\limp\;\stk\dR\cR\mo\psi(f,f_1,...,f_n)[x ]
\big)\,.
\dm
This implies $\qD\,x \:(\stk\dR\cR\mo\psi(f,f_1,...,f_n)[x ]).$
Then, by the inductive assumption, 
$
\stk\adR\aR\mo\psi(\kla f,\kla{f_1},...,\kla{f_n}),
$ 
hence  
$\stk\adR\aR\mo\vpi(\kla{f_1},...,\kla{f_n}),$ as required.\vom

\epF{Lemma}

To accomplish the proof of Theorem~\ref e, consider a closed 
\dd\cL formula $\vpi(r_1,...,r_n)$ with parameters 
$r_1,...,r_n\in\dR.$ 
We have to prove the equivalence 
\dm
\stk\dR\cR\mo\vpi(r_1,...,r_n)\;\leqv \;
\stk\adR\aR\mo\vpi(\ar_1,...,\ar_n)\,.
\dm
Let $f_i=c_{r_i},$ thus, $f_i\in F$ and $f_i(x )=r_i\zd\kaz x.$ 
Obviously $\vpi(f_1,...,f_n)[x ]$ coincides with 
$\vpi(r_1,...,r_n)$ for any $x \in\dN^A,$ hence 
$\vpi(r_1,...,r_n)$ is equivalent to 
$\qD\,x \:\vpi(f_1,...,f_n)[x ].$ 
On the other hand, by definition, $\ar_i=\kla{f_i}.$ 
Now the result follows by Lemma~\ref l.
\epf

\sekd{The iteration}
\label{s3}
                          
Theorem~\ref e yields a definable proper elementary extension 
$\stk\adR\aR$ of the structure $\stk\dR\cR.$ 
Yet this 
extension is not countably saturated due to the fact that 
the ultrapower $\adR$ was defined with maps concentrated on 
finite sets $u\sq A$ only. 
To fix this problem, we iterate the extension used above 
\dd\omi many times. 

Suppose that $\stk M\cM$ is an \dd\cL structure, so that 
$\cM$ consists of finitary relations on a set $M,$ and for any 
$E\in\cR$ there is a relation $E^\cM\in\cM$ of the same arity, 
associated with $E.$ 
Let $F_M$ be the set of all maps $f:\dN^A\to  M$ concentrated on 
finite sets $u\sq A.$ 
The structure $F_M/D=\stk\aM\acM,$ defined as in Section~\ref{s2}, 
but with 
the modified $F,$ will be called {\it the \dd Dultrapower\/} 
of $\stk M\cM.$ 
Theorem~\ref e remains true in this general setting: 
the map $x\mto \ax\;\,(x\in M)$  
is an elementary embedding of $\stk M\cM$ in $\stk\aM\acM$.

We define a sequence of \dd\cL structures 
$\stk{M_\al}{\cM_\al}\zt\al\le\omi,$ together with a system of 
elementary embeddings 
$e_{\al\ba}:\stk{M_\al}{\cM_\al}\to\stk{M_\ba}{\cM_\ba}\zt
\al<\ba\le\omi,$ so that
\ben
\enrm
\itla{x1}\msur
$\stk{M_0}{\cM_0}=\stk\dR\cR$;

\itla{x2}\msur
$\stk{M_{\al+1}}{\cM_{\al+1}}$ is the \dd Dultrapower of
$\stk{M_\al}{\cM_\al},$  that is, 
$\stk{M_{\al+1}}{\cM_{\al+1}}=F_\al/D,$ where 
$F_\al=F_{M_\al}$ consists of all functions 
$f:\dN^A\to M_\al$ concentrated on finite sets $u\sq A.$
In addition,
$e_{\al,\al+1}$ is the associated \dd{^\ast}embedding 
$\stk{M_\al}{\cM_\al}\to\stk{M_{\al+1}}{\cM_{\al+1}},$ 
while $e_{\ga,\al+1}=e_{\al,\al+1}\circ e_{\ga\al}$ for any 
$\ga<\al$ 
(in other words, 
$e_{\ga,\al+1}(x)= e_{\al,\al+1}(e_{\ga\al}(x))$ for all 
$x\in M_\al$);

\itla{x3}
if $\la\le\omi$ is a limit ordinal then $\stk{M_\la}{\cM_\la}$ 
is the direct limit of the structures 
$\stk{M_\al}{\cM_\al}\zt\al<\la.$ 
This can be achieved by the following steps: 
%
\ben
\enar
\itla{z1}\msur
$M_\la$ is defined  
as the set of all pairs $\ang{\al,x}$ such that $x\in M_\al$ 
and $x\nin\ran e_{\ga\al}$ for all $\ga<\al$.
  
\itla{z2}
If $E\in\cR$ is an \dd nary relation symbol then we define 
an \dd nary relation $E_\la$ on $M_\la$ as follows.  
Suppose that 
$\fx_i=\ang{\al_i,x_i}\in M_\la$ for $i=1,...,n.$  
Let $\al=\tsup{\ans{\al_1,...,\al_n}}$ and 
$z_i=e_{\al_i,\al}(x_i)$ for every $i,$ so that 
$\al_i\le\al<\la$ and $z_i\in M_\al.$ 
(Note that if $\al_i=\al$ then $e_{\al_i,\al}$ is the identity.)
Define $E_\la(\fx_1,...,\fx_n)$ iff 
$\stk{M_\al}{\cM_\al}\mo E(z_1,...,z_n).$  
  
\itla{z3}
Put $\cM_\la=\ans{E_\la:E\in\cR}$ -- 
then $\stk{M_\la}{\cM_\la}$ is an \dd\cL structure. 
  
\itla{z4}
Define an embedding $e_{\al\la}:M_\al\to M_\la\;\,(\al<\la)$ 
as follows. 
Consider any $x\in M_\al.$ 
If there is a least $\ga<\al$ such that there exists an element 
$y\in M_\ga$ with $x=e_{\ga\al}(y)$ then let 
$e_{\al\la}(x)=\ang{\ga,y}.$ 
Otherwise put $e_{\al\la}(x)=\ang{\al,x}$.
\een
\een

A routine verification of the following is left to the 
reader.

\bpro
\label p
If\/ $\al<\ba\le\omi$ then\/ $e_{\al\ba}$ is an elementary 
embedding of\/ $\stk{M_\al}{\cM_\al}$ to\/ 
$\stk{M_\ba}{\cM_\ba}$.\qed
\epro

Note that the construction of the sequence of models 
$\stk{M_\al}{\cM_\al}$ is definable, hence, so is the last 
member $\stk{M_{\omi}}{\cM_{\omi}}$ of the sequence. 
It remains to prove that the \dd\cL structure 
$\stk{M_{\omi}}{\cM_{\omi}}$ is countably saturated. 

This is also a simple argument. 
Suppose that, for any $k,$  $\vpi_k(p_k,x)$ is an 
\dd\cL formula with a single parameter $p_k\in M_{\omi}$ 
(the case of many parameters does not essentially differ 
from the case of one parameter), 
and there exists an element $x_k\in\ M_{\omi}$ such that 
$\bigwedge_{i\le k} \vpi_i(p_i,x_k)$ is true in 
$\stk{M_{\omi}}{\cM_{\omi}}$ --- in other words, we have 
$\stk{M_{\omi}}{\cM_{\omi}}\mo\vpi_i(p_i,x_k)$ whenever 
$k\ge i.$
Fix an ordinal $\ga<\omi$ such that for any $k\zi i$ there exist  
(then obviously unique) 
$y_k\zi q_i\in M_\ga$ with $x_k=e_{\ga\omi}(y_k)$ and 
$p_i=e_{\ga\omi}(q_i).$ 
Then $\vpi_i(q_i,y_k)$ is true in $\stk{M_\ga}{\cM_\ga}$ 
whenever $k\ge i$.

Fix $a\in A$ such that $D_a$ is a non-principal ultrafilter, 
that is, all cofinite subsets of $\dN$ belong to $D_a.$ 
Consider the structure $\stk{M_{\ga+1}}{\cM_{\ga+1}}$  as the 
\dd Dultrapower of $\stk{M_\ga}{\cM_\ga}.$ 
The corresponding set $F_\ga$ consists of all functions 
$f:\dN^A\to M_\ga$ concentrated on finite sets $u\sq A.$ 
In particular, the map $f(x )=y_k$ whenewer $x(a)=k$ 
belongs to $F_\ga.$ 
As any set of the form $\ans{k:k\ge i}$ belongs to $D_a,$ 
we have $\qD_a\, k\;(\stk{M_\ga}{\cM_\ga}\mo\vpi_i(q_i,y_k)),$ 
that is, 
$\qD\,x \in \dN^A\:(\stk{M_\ga}{\cM_\ga}\mo\vpi_i(q_i,f)[x ]),$
for any $i\in\dN.$ 
It follows, by Lemma~\ref l, that $\vpi_i(\aq_i,\fy)$ holds in 
$\stk{M_{\ga+1}}{\cM_{\ga+1}}$ for any $i,$ where 
$\aq_i=e_{\ga,\ga+1}(q_i)\in M_{\ga+1}$ while 
$\fy=\kla f\in M_{\ga+1}$ is the \dd Dequivalence class of $f$ 
in $F_\ga.$ 
Put $\fx=e_{\ga+1,\omi}(\fy);$ then 
$\vpi_i(p_i,\fx)$ is true in $\stk{M_{\omi}}{\cM_{\omi}}$
for any $i$ 
because obviously $p_i=e_{\ga+1,\omi}(\aq_i)\zd\kaz i$.\vom

\qeD{Theorem~\ref t}

\sekd{Varia}
\label{s4}
                          
By appropriate modifications of the constructions, the following 
can be achieved:

\ben
\item
For any given infinite cardinal $\ka,$ a \dd\ka saturated 
elementary extension of $\dR,$ definable with $\ka$ 
as the only parameter of definition.

\item
A {\it special\/} elementary extension of $\dR,$ 
of as large cardinality as desired. 
For instance, take, in stage $\al$ of the construction 
considered in Section~\ref{s3}, ultrafilters on 
$\beth_{\al}.$ 
Then the result will be a definable special structure of 
cardinality $\beth_{\omi}.$ 
Recall that special models of equal cardinality are 
isomorphic \cite[Theorem~5.1.17]{ck}.   
Therefore, such a modification admits an explicit 
model-theoretical characterization up to isomorphism.

\item
A class-size definable elementary extension of $\dR,$ 
\dd\ka saturated for any cardinal $\ka$.

\item
A class-size definable elementary extension of 
the whole set universe, 
\dd\ka saturated for any cardinal $\ka.$ 
(Note that this cannot be strengthened to 
\dd\Ord saturation, \ie, saturation with respect to all  
class-size families. 
For instance, \dd{\Ord^M}saturated elementary extensions of a 
minimal transitive model $M\mo\ZFC,$ definable in $M,$ 
do not exist --- see \cite[Theorem~2.8]{hyp}.)
\een

The authors thank the anonimous referee for valuable 
comments and corrections.

\renek{\refname}
\vspace{10mm}

\noi
Vladimir Kanovei\\
Institute for information transmission problems (IPPI),  
Russian academy of sciences, 
Bol. Karetnyj Per. 19, 
Moscow 127994, Russia\\
E-mail address : {\tt kanovei@mccme.ru}\\ 

\noi
Saharon Shelah\\
Institute of Mathematics, The Hebrew University of Jerusalem, 91904 
Jerusalem, Israel, and Department of Mathematics, Rutgers University, 
New Brunswick, NJ 08854, USA\\ 
E-mail address : {\tt shelah@math.huji.ac.il}\\ 
URL: \verb" http://www.math.rutgers.edu/~shelah"

\end{document}